\documentclass[11pt]{article}
\usepackage{amsmath, amsfonts, amsthm, amssymb, color}
\usepackage{graphicx}
\usepackage{float}
\usepackage{verbatim}

\allowdisplaybreaks

\numberwithin{equation}{section}

\newtheorem{lemma}{Lemma}[section]

\newtheorem{thm}[lemma]{Theorem}

\theoremstyle{definition}

\newtheorem{conjecture}[lemma]{Conjecture}

\theoremstyle{remark}
\newtheorem{remark}[lemma]{Remark}

\def\R{\mathbb{R}}
\def\N{\mathbb{N}}

\def\a{\mathbf{a}}
\def\b{\mathbf{b}}

\def\I{\mathcal{I}}

\numberwithin{equation}{section} \numberwithin{table}{section}

\title{Iterated function systems with super-exponentially close cylinders}
\author{Simon Baker\\ \\
\emph{School of Mathematics,} \\ \emph{University of Birmingham,} \\ \emph{Birmingham,  B15 2TT, UK.} \\ Email: simonbaker412@gmail.com\\}

\date{\today}
\begin{document}
\maketitle

\begin{abstract}
Several important conjectures in Fractal Geometry can be summarised as follows: If the dimension of a self-similar measure in $\R$ does not equal its expected value, then the underlying iterated function system contains an exact overlap. In recent years significant progress has been made towards these conjectures. Hochman proved that if the Hausdorff dimension of a self-similar measure in $\R$ does not equal its expected value, then there are cylinders which are super-exponentially close at all small scales. Several years later, Shmerkin proved an analogous statement for the $L^q$ dimension of self-similar measures in $\R$. With these statements in mind, it is natural to wonder whether there exist iterated function systems that do not contain exact overlaps, yet there are cylinders which are super-exponentially close at all small scales.
In this paper we show that such iterated function systems do exist. In fact we prove much more. We prove that for any sequence $(\epsilon_n)_{n=1}^{\infty}$ of positive real numbers, there exists an iterated function system $\{\phi_i\}_{i\in \I}$ that does not contain exact overlaps and  $$\min\left\{|\phi_{\a}(0)-\phi_{\b}(0)|: \a,\b\in \I^n,\, \a\neq \b,\, r_{\a}=r_{\b}\right\}\leq \epsilon_n$$ for all $n\in \N.$\\



\end{abstract}

\section{Introduction}
Let $\Phi=\{\phi_i(x)=r_ix+a_i\}_{i\in \I}$ be a finite set of contracting similarities on $\R$. We call $\Phi$ an iterated function system or IFS for short. An important result due to Hutchinson \cite{Hut} states that for any IFS, there exists a unique non-empty compact set $X\subseteq \mathbb{R}$ such that $$X=\bigcup_{i\in \I}\phi_i(X).$$ We call $X$ the self-similar set of $\Phi$. One of the most important problems within Fractal Geometry is determining the metric properties of self-similar sets. In particular, we are often interested in determining their Hausdorff dimension. The following upper bound is known to hold for any IFS: Let $\dim_{S}\Phi$ be the unique solution to $\sum_{i\in \I}|r_i|^s=1,$ then 
\begin{equation}
\label{expected set dimension}
\dim_{H}X\leq \min\{\dim_{S}\Phi,1\}.
\end{equation}It is believed that equality in \eqref{expected set dimension} holds generically. Determining what conditions prevent equality is an active area of research. To make progress with this problem one studies self-similar measures. Let $\Phi$ be an IFS and $(p_i)_{i\in \I}$ be a probability vector, then there is a unique Borel probability measure $\mu$ supported on $X$ satisfying $$\mu=\sum_{i\in\I} p_i\cdot  \mu\circ \phi_i^{-1}.$$ The measure $\mu$ is called the self-similar measure associated to $\Phi$ and $(p_i).$ For a Borel probability measure $\nu$ on $\R$, we define the dimension of $\nu$ to be $$\dim \nu:=\inf\{\dim_{H}A:\nu(A)=1\}.$$ There are several other notions of dimension for measures, but for self-similar measures all of the main ones coincide. This is a consequence of the exact dimensionality of self-similar measures \cite{FengHu}. The following upper bound for the dimension of self-similar measures is well known:
\begin{equation}
\label{expected measure dimension}
\dim \mu\leq \min\left\{\frac{\sum_{i\in \I}p_i\log p_i}{\sum_{i\in \I}p_i\log |r_i|},1\right\}.
\end{equation}
As is the case for \eqref{expected set dimension}, it is believed that equality in \eqref{expected measure dimension} holds generically. Note that taking our probability vector to equal $(|r_i|^{\dim_{S}(\Phi)})_{i\in \I},$ it can be shown that equality in \eqref{expected measure dimension} implies equality in \eqref{expected set dimension}.

A self-similar measure may distribute mass within $X$ in an non-uniform manner. Understanding the irregularities in the distribution of a measure is the topic of multifractal analysis. One of the standard tools for describing the multifractal properties of a measure is the $L^q$ dimension. Let $q>1$ and $\nu$ be a Borel probability measure on $\mathbb{R},$ we define the $L^q$ dimension of $\nu$ to be
$$D(\nu,q):=\liminf_{n\to\infty}\frac{-\log \sum_{j\in \mathbb{Z}}\nu([j\cdot2^{-n},(j+1)\cdot2^{-n}))}{(q-1)n}.$$ The following upper bound for the $L^q$ dimension of a self-similar measure always holds:
\begin{equation}
\label{Lq upper bound}
D(\mu,q)\leq \min\left\{\frac{T(\mu,q)}{q-1},1\right\}.
\end{equation}Here $T(\mu,q)$ is the unique solution to the equation $$\sum_{i\in \I}p_i^q|r_i|^{-T(\mu,q)}=1.$$ Equality in \eqref{Lq upper bound} is believed to hold generically.

One of the most important and active areas of Fractal Geometry is concerned with determining when we have equality in \eqref{expected set dimension}, \eqref{expected measure dimension}, and \eqref{Lq upper bound}, see \cite{Hochman2, Hochman, Shm, Shm2} and the references therein for more on this topic. One can construct examples where we have a strict inequality in \eqref{expected set dimension}, \eqref{expected measure dimension}, and \eqref{Lq upper bound}, by building an IFS in such a way that it contains an exact overlap. Recall that we say that an IFS contains an exact overlap if there exists $\a,\b\in \I^n$ such that $\a\neq \b$ and $\phi_{\a}=\phi_{\b}.$ Here and throughout, for $\a=(a_1,\ldots,a_n)$ we let $\phi_{\a}$ denote $\phi_{a_1}\circ \cdots \circ \phi_{a_n}$ and $r_{\a}$ denote the scaling ratio of $\phi_{\a}$. If an IFS contains an exact overlap then we can remove a similarity from $\{\phi_{\a}\}_{\a\in \I^n}$ without changing the self-similar set. This observation can result in a strict inequality in \eqref{expected set dimension}, \eqref{expected measure dimension}, and \eqref{Lq upper bound}. An important folklore conjecture in this area states that this is the only way a strict inequality can occur. 

\begin{conjecture}
\label{folklore}
Let $\Phi$ be an IFS and $(p_i)_{i\in I}$ be a probability vector. If we have strict inequality in one of \eqref{expected set dimension}, \eqref{expected measure dimension}, or \eqref{Lq upper bound}, then our IFS contains an exact overlap.
\end{conjecture}

Significant progress on Conjecture \ref{folklore} has been made in recent years. Hochman in \cite{Hochman} and Shmerkin in \cite{Shm2} proved results to the effect that if we have a strict inequality in one of \eqref{expected set dimension}, \eqref{expected measure dimension}, or \eqref{Lq upper bound}, then our IFS comes extremely close to containing an exact overlap. To properly describe this notion of closeness we introduce the following distance function. Given $\a,\b\in \I^n, $ let $$d(\a,\b):=\begin{cases}
\infty,\, &r_{\a}\neq r_{\b}\\
|\phi_{\a}(0)-\phi_{\b}(0)|,\, &r_{\a}= r_{\b}.
\end{cases}$$ Importantly $d(\a,\b)=0$ if and only if $\phi_{\a}=\phi_{\b}.$ For any $n\in \N$ we let $$\Delta_n:=\min\{d(\a,\b):\, \a,\b\in \I^n,\, \a\neq \b\}.$$ The results of Hochman and Shmerkin can now be stated more accurately as follows. 
\begin{thm}
	\label{HocShm}Let $\Phi$ be an IFS.  
	\begin{itemize}
		\item \cite[Theorem 1.1]{Hochman} If $\limsup_{n\to\infty}\frac{\log \Delta_n}{n}>-\infty$ then we have equality in \eqref{expected set dimension} and \eqref{expected measure dimension} for all self-similar measures.
		\item \cite[Theorem 6.6]{Shm2} If $\limsup_{n\to\infty}\frac{\log \Delta_n}{n}>-\infty$ then we have equality in \eqref{Lq upper bound} for all self-similar measures and $q>1$.
	\end{itemize}
\end{thm}We mention for completeness that the first statement in Theorem \ref{HocShm} can be recovered from the second using the fact that for self-similar measures $\lim_{q\to 1^{+}}\frac{T(\mu,q)}{q-1}=\frac{\sum_{i\in \I}p_i\log p_i}{\sum_{i\in \I}p_i\log |r_i|}$. Theorem \ref{HocShm} tells us that if we have a strict inequality in one of \eqref{expected set dimension}, \eqref{expected measure dimension}, or \eqref{Lq upper bound}, then $\Delta_n\to 0$ super-exponentially fast. This fact can be used to provide a positive answer to Conjecture \ref{folklore} under the additional assumption that the scaling ratios and translation parameters are given by algebraic numbers, see Theorem 1.5 and Lemma 5.10 from \cite{Hochman}. The results summarised in Theorem \ref{HocShm} have subsequently been applied to many other problems from Fractal Geometry. Perhaps most notable amongst these applications are a proof of a conjecture of Furstenberg on projections of the one-dimensional Sierpinski gasket \cite{Hochman}, a proof that outside of a set of zero Hausdorff dimension the Bernoulli convolution is absolutely continuous with a density in $L^q$ for all finite $q>1$ \cite{Shm2}, and a proof that the Bernoulli convolution has dimension $1$ whenever the contraction parameter is transcendental \cite{Varju2}. 

For many parameterised families of iterated function systems it can be shown that the set of parameters for which $\Delta_n$ converges to zero super-exponentially fast is very small. One might hope that there exists a dichotomy which states that for an arbitrary IFS $\Phi$, either $\Phi$ contains an exact overlap or $\limsup_{n\to\infty}\frac{\log \Delta_n}{n}>-\infty.$ If such a dichotomy were to exist then Theorem \ref{HocShm} would imply Conjecture \ref{folklore}. In this paper we show that no such dichotomy exists. We give the first example of an IFS that does not contain exact overlaps yet $\lim_{n\to\infty}\frac{\log \Delta_n}{n}=-\infty$\footnote{After completing this paper the author was made aware of an upcoming article of B\'{a}r\'{a}ny and K\"{a}enm\"{a}ki \cite{BarKae}. In this article, using different methods, they also show that there exists an IFS without exact overlaps such that $\lim_{n\to\infty}\frac{\log \Delta_n}{n}=-\infty$.}. In fact we prove much more, we show that it is possible for an IFS to have no exact overlaps yet $\Delta_n$ converges to zero arbitrarily fast.

\begin{thm}
	\label{Main theorem}
Let $(\epsilon_n)_{n=1}^{\infty}$ be an arbitrary sequence of positive real numbers. Then there exists an IFS $\Phi$ without exact overlaps such that $\Delta_n\leq \epsilon_n$ for all $n\in \N$.
\end{thm}
The IFS we will construct in our proof of Theorem \ref{Main theorem} satisfies $\dim_{S}\Phi>1$. It is possible to adapt the construction in such a way that the IFS obtained satisfies $\dim_{S}\Phi<1.$ We will detail the necessary changes in Remark \ref{remark} after we have given our proof of Theorem \ref{Main theorem}. 
\section{Proof of Theorem \ref{Main theorem}}
\subsection{Continued fractions}
Before giving our proof of Theorem \ref{Main theorem} it is useful to recall some properties of continued fractions. Proofs of the properties stated below can be found in \cite{Bug} and \cite{Cas}.  

For any $s\in [0,1]\setminus \mathbb{Q},$ there exists a unique sequence $(\zeta_m)_{m=1}^{\infty}\in\mathbb{N}^{\mathbb{N}}$ such that 
$$ s=\cfrac{1}{\zeta_1+\cfrac{1}{\zeta_2 +\cfrac{1}
		{\zeta_3 + \cdots }}}.$$
We call the sequence $(\zeta_m)$ the continued fraction expansion of $s$. Suppose $s$ has continued fraction expansion $(\zeta_m),$ then for each $m\in\mathbb{N}$ we let
$$ \frac{p_m}{q_m}:=\cfrac{1}{\zeta_1+\cfrac{1}{\zeta_2 +\cfrac{1}
		{\zeta_3 + \cdots \cfrac{1}
			{\zeta_m }}}}.$$ We call $p_m/q_m$ the $m$-th partial quotient of $s$. The sequence of partial quotients of $s$ satisfy the following properties:
\begin{itemize}
	\item For any $m\in\mathbb{N}$ we have 
	\begin{equation}
	\label{property1}\frac{1}{q_m(q_{m+1}+q_m)}<\left|s-\frac{p_m}{q_m}\right|<\frac{1}{q_mq_{m+1}}.
	\end{equation}
	\item If we set $p_{-1}=1, q_{-1}=0, p_0=0, q_0=1$, then for any $m\geq 1$ we have 
	\begin{align}
	\label{property2}
	p_m&=\zeta_m p_{m-1}+p_{m-2}\\
	q_m&=\zeta_m q_{m-1}+q_{m-2}. \nonumber
	\end{align}
	\item If $q< q_{m+1}$ then \begin{equation}
	\label{property3}
	|qs-p|\geq |q_{m}s-p_m|
	\end{equation}for any $p\in\mathbb{Z}$.
\end{itemize}  
Given a finite word $(\zeta_1,\ldots,\zeta_n)\in \N^n,$ we let $$C_{\zeta_1,\ldots,\zeta_n}:=\left\{s\in [0,1]\setminus\mathbb{Q}:\textrm{ the continued fraction expansion of }s \textrm{ begins with }(\zeta_1,\ldots,\zeta_n)\right\}$$

\subsection{Our construction}
We are now in a position to prove Theorem \ref{Main theorem}. Our proof is based upon ideas from \cite[Section 5]{Bak} where the author studied a family of iterated function systems whose overlapping behaviour is determined by the Diophantine properties of an underlying translation parameter.

In the rest of this section $(\epsilon_n)_{n=1}^{\infty}$ is fixed. Without loss of generality we may assume $\epsilon_n\leq 8^{-n}$ for all $n\in\N$. The IFS we construct will be of the form 
\begin{align*}
\Phi_{s,t}=\Big\{&\phi_1(x)=\frac{x}{2},\,\phi_2(x)=\frac{x+1}{2},\, \phi_3(x)=\frac{x+s}{2},\, \phi_4(x)=\frac{x+1+s}{2},\, \phi_5(x)=\frac{x+t}{2},\\
&\phi_6(x)=\frac{x+1+t}{2}\Big\}.
\end{align*}Here $s,t$ are two real numbers contained in $[0,1]\setminus\mathbb{Q}$. We will define $(\zeta_m)$ and $(\zeta_m')$ in such a way that taking $s$ and $t$ to be the real numbers with continued fraction expansions $(\zeta_m)$ and $(\zeta_m')$ respectively, then $\Delta_n\leq \epsilon_n$ for all $n\in \N$ and $\Phi_{s,t}$ will not contain an exact overlap.

The following observations will be useful in our proof. For any $n\in \N$ we have 
\begin{align}
\label{zero images}
&\left\{\phi_{\a}(0):\a\in\{1,\ldots,6\}^n\right\}\nonumber\\
=&\left\{\sum_{i=1}^n\frac{\omega_i}{2^i}+s\sum_{i=1}^n\frac{\delta_i}{2^i}+t\sum_{i=1}^n\frac{\delta_i'}{2^i}:(\omega_i),(\delta_i),(\delta_i')\in\{0,1\}^n,\, \delta_i\cdot \delta_i'=0\, \forall 1\leq i\leq n \right\}.
\end{align}Equation \eqref{zero images} is easily seen to hold when $n=1$. The general case can then be deduced by an induction argument. The following inclusions, which are a consequence of \eqref{zero images}, will be very useful to us. For any $n\in \N$ we have 
\begin{equation*}
\left\{\sum_{i=1}^n\frac{\omega_i}{2^i}+s\sum_{i=1}^n\frac{\delta_i}{2^i}:(\omega_i),(\delta_i)\in\{0,1\}^n\right\}\subseteq \left\{\phi_{\a}(0):\a\in\{1,\ldots,6\}^n\right\}
\end{equation*}and 
\begin{equation*}
\left\{\sum_{i=1}^n\frac{\omega_i}{2^i}+t\sum_{i=1}^n\frac{\delta_i'}{2^i}:(\omega_i),(\delta_i')\in\{0,1\}^n\right\}\subseteq \left\{\phi_{\a}(0):\a\in\{1,\ldots,6\}^n\right\}.
\end{equation*} Equivalently one has 
\begin{equation}
\label{s inclusion}
\left\{\frac{1}{2^n}(p+sq):0\leq p,q\leq 2^{n}-1\right\}\subseteq \left\{\phi_{\a}(0):\a\in\{1,\ldots,6\}^n\right\}
\end{equation} and 
\begin{equation}
\label{t inclusion}
\left\{\frac{1}{2^n}(p+tq):0\leq p,q\leq 2^{n}-1\right\}\subseteq \left\{\phi_{\a}(0):\a\in\{1,\ldots,6\}^n\right\}.
\end{equation}The following lemma is an immediate consequence of \eqref{s inclusion}, \eqref{t inclusion}, and the irrationality of $s$ and $t$.

\begin{lemma}
	\label{Delta bound}
For any $n\in \N$ we have
	$$	\Delta_n\leq \min\left\{ \min_{\stackrel{0\leq p,q\leq 2^n-1}{(p,q)\neq (0,0)}}\{\frac{1}{2^n}|qs-p|\},\min_{\stackrel{0\leq p,q\leq 2^n-1}{(p,q)\neq (0,0)}}\{\frac{1}{2^n}|qt-p|\}\right\}.$$
\end{lemma}
Before constructing $(\zeta_m)$ and $(\zeta_m')$ we say a few words on our proof. Suppose $(\zeta_1,\ldots,\zeta_k)$ has been defined and we let $L_k$ be such that $2^{L_k-1}-1\leq q_k< 2^{L_k}-1$. If we choose $\zeta_{k+1}$ to be some extraordinarily large integer, then \eqref{property1} and \eqref{property2} imply that $|q_ks-p_k|$ is extremely small. Moreover, if we also chose $\zeta_{k+1}$ in a way that depends upon $(\epsilon_n),$ then Lemma \ref{Delta bound} can be used to prove that $\Delta_n\leq \epsilon_n$ for all $L_k\leq n\leq L_k+K_k$ for some $K_k\in \N$. Unfortunately it is not the case that $2^{L_k+K_k}>q_{k+1}.$ If we were to repeat this process of picking large digits in the continued fraction expansion of $s$, then we would obtain infinitely many disjoint closed intervals in $[1,\infty),$ such that if $n$ is in one of these intervals we have $\Delta_n\leq \epsilon_n$. The key observation here is that all of the above remains true if $s$ is replaced by $t$. We can choose $(\zeta_m')$ in such a way that we obtain infinitely many disjoint closed intervals in $[1,\infty),$ such that if $n$ is in one of these intervals we have $\Delta_n\leq \epsilon_n$. In our proof of Theorem \ref{Main theorem} we define $(\zeta_m)$ and $(\zeta_m')$ in such a way that the union of the intervals corresponding to $s$ and $t$ covers $[1,\infty)$. This ensures $\Delta_n\leq \epsilon_n$ for all $n\in \N$. If $\Phi_{s,t}$ contains an exact overlap then $t$ is equal to the image of $s$ under some affine map with rational coefficients. This means that good rational approximations for $s$ can be translated into good rational approximations for $t$ at approximately the same scale. We will define $s$ and $t$ in such a way that they have extremely good rational approximations at very different scales. This property and the remark stated previously will ensure that exact overlaps are not possible. \\

\noindent \textbf{Initial step.}
We start our construction by letting $\zeta_1=1$ and choosing $N\in\mathbb{N}$ to be some large number. By \eqref{property2} we know that $p_1=1$ and $q_1=1$. We let $$\zeta_2=\left\lceil\frac{1}{\min_{1\leq n\leq N}\epsilon_n}\right\rceil.$$ Then by \eqref{property1} and \eqref{property2}, we see that for any $s\in C_{\zeta_1,\zeta_2}$ we have
\begin{equation*}
|s-1|=|q_1s-p_1|<\frac{1}{q_{2}}\leq \frac{1}{\zeta_2}\leq \min_{1\leq n\leq N}\epsilon_n.
\end{equation*} Which by Lemma \ref{Delta bound} implies $$\Delta_n\leq \epsilon_n$$ for all $1\leq n\leq N$ whenever $s\in C_{\zeta_1,\zeta_2}$

Let $L_2\in\N$ be such that $$2^{L_2-1}-1\leq q_2<2^{L_2}-1.$$ We now choose $\zeta_1'=1$ and let
$$\zeta_{2}'=\left\lceil\frac{1}{\min_{1\leq n\leq L_2}\epsilon_n}\right\rceil.$$  Then $q_1'=1$ and $p_1'=1$. By \eqref{property1} and \eqref{property2} it follows that for any $t\in C_{\zeta_1',\zeta_2'}$ we have
\begin{equation*}
|t-1|=|q_1't-p_1'|<\frac{1}{q_{2}'}\leq \frac{1}{\zeta_2'}\leq \min_{1\leq n\leq L_2}\epsilon_n.
\end{equation*} Which by Lemma \ref{Delta bound} implies that $$\Delta_n\leq \epsilon_n$$ for all $1\leq n\leq L_2$ whenever $t\in C_{\zeta_1,\zeta_2}.$ 

If we let $M_2\in\N$ be such that $$2^{M_{2}-1}-1\leq q_{2}'<2^{M_{2}}-1,$$ then we have $2L_2<M_2$. This follows from the definition of $M_2$ and because $$q_{2}'\geq \zeta_{2}'\geq \frac{1}{\epsilon_{L_2}}\geq 8^{L_2}.$$ Here we used our assumption $\epsilon_n\leq 8^{-n}$ for all $n\in \N$.  \\

\noindent \textbf{Iterative step.} Suppose we have constructed $(\zeta_m)_{m=1}^{k+1}$ and $(\zeta_m')_{m=1}^{k+1}$. For each $2\leq m\leq k+1$ let $L_m\in \mathbb{N}$ be such that  $$2^{L_{m}-1}-1\leq q_{m}<2^{L_{m}}-1,$$ and  $M_m\in \N$ be such that $$2^{M_{m}-1}-1\leq q_{m}'<2^{M_{m}}-1.$$ Assume that the following properties are satisfied:
\begin{enumerate}
	\item For any $s\in C_{\zeta_1,\ldots,\zeta_{k+1}}$ and $t\in C_{\zeta_1',\ldots,\zeta_{k+1}'}$ we have $\Delta_n\leq \epsilon_n$ for all $1\leq n\leq L_{k+1}.$
	\item For each $2\leq m\leq k+1$ we have $2L_m<M_m.$
\end{enumerate}	By our initial step we know that these properties are satisfied when $k=1$. We now define $\zeta_{k+2}$ and $\zeta_{k+2}'$ in such a way that these properties are still satisfied by $(\zeta_m)_{m=1}^{k+2}$ and $(\zeta_m')_{m=1}^{k+2}$.\\

Let 
\begin{equation*}
\zeta_{k+2}=\left\lceil \frac{1}{\min_{L_{k+1}\leq n\leq M_{k+1}}\epsilon_n}\right\rceil.
\end{equation*} Then by \eqref{property1} and \eqref{property2} \begin{equation}
\label{refer to}
|q_{k+1}s-p_{k+1}|\leq \frac{1}{\zeta_{k+2}}\leq  \min_{L_{k+1}\leq n\leq M_{k+1}}\epsilon_n.
\end{equation} Which by Lemma \ref{Delta bound} implies 
\begin{equation}
\label{s working}
\Delta_n\leq \epsilon_n \textrm{ for all } L_{k+1}\leq n \leq M_{k+1}
\end{equation}
 when $s\in C_{\zeta_1,\ldots,\zeta_{k+2}}.$ Define $L_{k+2}\in \N$ to be the solution to $$2^{L_{k+2}-1}-1\leq q_{k+2}<2^{L_{k+2}}-1.$$ Let $$\zeta_{k+2}'=\left\lceil\frac{1}{\min_{M_{k+1}\leq n\leq L_{k+2}}\epsilon_n}\right\rceil.$$ Then by \eqref{property1} and \eqref{property2} $$|q_{k+1}'t-p_{k+1}'|\leq \frac{1}{\zeta_{k+2}'}\leq \min_{M_{k+1}\leq n\leq L_{k+2}}\epsilon_n.$$ Which by Lemma \ref{Delta bound} implies
 \begin{equation}
 \label{t working}
 \Delta_n\leq \epsilon_n \textrm{ for all } M_{k+1}\leq n \leq L_{k+2}
 \end{equation} when $t\in C_{\zeta_1',\ldots,\zeta_{k+2}'}.$ Combining \eqref{s working} and \eqref{t working} with our underlying assumptions, we see that  
 \begin{equation}
 \label{winner}
 \Delta_n\leq \epsilon_n \textrm{ for all } 1\leq n\leq L_{k+2} 
 \end{equation}when $s\in C_{\zeta_1,\ldots,\zeta_{k+2}}$ and $t\in C_{\zeta_1',\ldots,\zeta_{k+2}'}.$ 
 
 If $M_{k+2}\in \N$ is the unique solution to $$2^{M_{k+2}-1}-1\leq q_{k+2}'<2^{M_{k+2}-1}-1,$$ then we have $2L_{k+2}<M_{k+2}.$ As in our initial step this is a consequence of our underlying assumption $\epsilon_n\leq 8^{-n}$ for all $n\in \N$. We have therefore defined $\zeta_{k+2}$ and $\zeta_{k+2}'$ in such a way that properties $1.$ and $2$. are satisfied by $(\zeta_m)_{m=1}^{k+2}$ and $(\zeta_m')_{m=1}^{k+2}$.\\
 
Clearly our iterative step can be repeated indefinitely, and doing so yields two sequences $(\zeta_m)$ and $(\zeta_m').$ The parameter $L_m$ tends to infinity. Therefore it follows from \eqref{winner} that if $s$ and $t$ have continued fraction expansions $(\zeta_m)$ and $(\zeta_m')$ respectively, then $$\Delta_n\leq \epsilon_n \textrm{ for all }n\in \N.$$ To complete our proof of Theorem \ref{Main theorem}, it remains to show that for the $s$ and $t$ constructed above the IFS $\Phi_{s,t}$ does not contain an exact overlap. To do this we will use the property $2L_m<M_m$ for all $m\geq 2$. \\

\noindent \textbf{No exact overlaps.}
Let $s$ and $t$ be as above. Assume $\Phi_{s,t}$ contains an exact overlap. By \eqref{zero images} there exists $(\omega_{i,1}),(\delta_{i,1}),(\delta_{i,1}'),(\omega_{i,2}),(\delta_{i,2}),(\delta_{i,2}')\in\{0,1\}^n$ such that, either $(\omega_{i,1})\neq (\omega_{i,2})$, $(\delta_{i,1})\neq (\delta_{i,2}),$ or  $(\delta_{i,1}')\neq (\delta_{i,2}'),$ and \begin{equation}
\label{exact overlap} \sum_{i=1}^n\frac{\omega_{i,1}}{2^i}+s\sum_{i=1}^n\frac{\delta_{i,1}}{2^i}+t\sum_{i=1}^n\frac{\delta_{i,1}'}{2^i}=\sum_{i=1}^n\frac{\omega_{i,2}}{2^i}+s\sum_{i=1}^n\frac{\delta_{i,2}}{2^i}+t\sum_{i=1}^n\frac{\delta_{i,2}'}{2^i}.
\end{equation}It follows from the irrationality of $s$ and $t$ that if \eqref{exact overlap} holds, then there must exist $a,b,c,d\in \mathbb{Z}$ such that $a\neq 0, b\neq 0, d\neq 0,$ and \begin{equation}
\label{rational relation}
t=\frac{a}{b}s+\frac{c}{d}.
\end{equation} By \eqref{property1} and \eqref{property3} we know that for any $m\geq 2$ we have $$\left|t-\frac{p}{q}\right|>\frac{1}{2q_{m}'^{2}}$$ for any $(p,q)\in \mathbb{Z}\times \N$ such that $q<q_{m}'.$ This in turn implies that for any $m\geq 2$ we have 
\begin{equation}
\label{step1}
\left|t-\frac{p}{q}\right|>\frac{1}{2^{2M_{m}+1}}
\end{equation} for any $(p,q)\in \mathbb{Z}\times \N$ such that $q<q_{m}'.$ By \eqref{refer to} and our assumption $\epsilon_n\leq 8^{-n},$ we know that $$\left|s-\frac{p_{m}}{q_{m}}\right|\leq 8^{-M_{m}}$$ for all $m\geq 2$. Therefore $$\left|\frac{a}{b}s+\frac{c}{d}-\left(\frac{ap_{m}}{bq_{m}}+\frac{c}{d}\right)\right|\leq \frac{a\cdot 8^{-M_{m}}}{b}.$$ Which by \eqref{rational relation} implies 
\begin{equation}
\label{step2}
\left|t-\frac{ap_{m}d+cbq_{m}}{bq_{m}d}\right|\leq \frac{a\cdot 8^{-M_{m}}}{b}.
\end{equation}Recall that $q_{m}\leq 2^{L_{m}}$ and $q_{m}'\geq 2^{M_{m}-4}$ for $m\geq 2$. Since $2L_m<M_m$ for $m\geq 2$ and $L_m$ tends to infinity, it follows that $bq_{m}d<q_{m}'$ for all $m$ sufficiently large. Therefore \eqref{step2} contradicts \eqref{step1} for $m$ sufficiently large. Thus we have obtained our contradiction and may conclude that $\Phi_{s,t}$ has no exact overlaps.

\begin{remark}
	\label{remark}
The IFS constructed in our proof of Theorem \ref{Main theorem} satisfies $\dim_{S}\Phi>1$. We now explain how the construction can be adapted to ensure $\dim_{S}\Phi<1.$ 

Consider the family of IFSs given by \begin{align*}\tilde{\Phi}_{s,t}=\Big\{&\phi_1(x)=\frac{x}{8},\,\phi_2(x)=\frac{x+1}{8},\, \phi_3(x)=\frac{x+s}{8},\, \phi_4(x)=\frac{x+1+s}{8},\, \phi_5(x)=\frac{x+t}{8},\\
&\phi_6(x)=\frac{x+1+t}{8}\Big\}.
\end{align*} Where once again $s,t\in [0,1]\setminus \mathbb{Q}$. For this family the appropriate analogue of Lemma \ref{Delta bound} is the following.

\begin{lemma}
	\label{Delta bound2}
	For any $n\in \N$ we have
	\begin{align*}
	\Delta_n\leq \min\Bigg\{& \min_{\stackrel{(\omega_i),(\delta_i)\in\{0,1\}^n}{\text{Either } (\omega_i)\neq 0^n \text{ or }(\delta_i)\neq 0^n}}\left\{\frac{1}{8^n}\left|s\sum_{i=1}^n\delta_i8^{n-i}-\sum_{i=1}^n\omega_i8^{n-i}\right|\right\}\\
	&,\min_{\stackrel{(\omega_i),(\delta_i')\in\{0,1\}^n}{\text{Either } (\omega_i)\neq 0^n \text{ or }(\delta_i')\neq 0^n}}\left\{\frac{1}{8^n}\left|t\sum_{i=1}^n\delta_i'8^{n-i}-\sum_{i=1}^n\omega_i8^{n-i}\right|\right\}\Bigg\}.
	\end{align*}
\end{lemma}
The important difference between Lemma \ref{Delta bound} and Lemma \ref{Delta bound2}, is that in Lemma \ref{Delta bound2}, the integers we may use to bound $\Delta_n$ must have base $8$ expansion consisting solely of zeros and ones. The reason our original proof doesn't immediately work for the family $\tilde{\Phi}_{s,t}$ is because it is not necessarily the case that the $p_m,q_m,p_m'$ and $q_m'$ constructed have base $8$ expansion consisting solely of zeros and ones. This issue can be overcome if we set $\zeta_1=\zeta_1'=1$ as in our original proof, and insist that for $m\geq 2$ the $\zeta_m$ and $\zeta_m'$ constructed are also suitably large powers of $8$. The recurrence relations $p_m=\zeta_m p_{m-1}+p_{m-2}$ and $q_m=\zeta_m q_{m-1}+q_{m-2},$ together with the initial conditions $p_{-1}=1, q_{-1}=0, p_0=0, q_0=1,$ can then be combined to ensure that the $p_m,q_m,p_m'$ and $q_m'$ constructed always have base $8$ expansion consisting solely of zeros and ones. With this issue overcome, the proof is now essentially identical. 
\end{remark}


\begin{remark}
If we assume $\epsilon_n\to 0$ super-exponentially fast, then Lemma 5.10 from \cite{Hochman} implies that one of the $s$ and $t$ constructed in Theorem \ref{Main theorem} is transcendental. A close examination of the proof of Theorem \ref{Main theorem} shows that for such an $(\epsilon_n)$ the $s$ and $t$ constructed are Liouville numbers. It is a well known fact that all Liouville numbers are transcendental.  
\end{remark}

\begin{remark}
Instead of using the family of IFSs $\{\Phi_{s,t}\}_{s,t}$ in our proof of Theorem \ref{Main theorem}, we could have instead used the family $\{\Phi_{s,t}'\}_{s,t}$, where 
\begin{align*}
\Phi_{s,t}'=\Big\{&\phi_1(x)=\frac{x}{2},\,\phi_2(x)=\frac{x+1}{2},\, \phi_3(x)=\frac{x+s}{2},\, \phi_4(x)=\frac{x+1+s}{2},\, \phi_5(x)=\frac{x+t}{2},\\
&\phi_6(x)=\frac{x+1+t}{2},\, \phi_7(x)=\frac{x+s+t}{2},\, \phi_8(x)=\frac{x+1+s+t}{2}\Big\}.
\end{align*}We chose $\{\Phi_{s,t}\}_{s,t}$ in our proof as we wanted to minimise the numbers of similarities appearing in our example. For $\Phi_{s,t}'$ we have 
$$\left\{\phi_{\a}(0):\a\in\{1,\ldots,8\}^n\right\}
=\left\{\sum_{i=1}^n\frac{\omega_i}{2^i}+s\sum_{i=1}^n\frac{\delta_i}{2^i}+t\sum_{i=1}^n\frac{\delta_i'}{2^i}:(\omega_i),(\delta_i),(\delta_i')\in\{0,1\}^n\right\}.$$
The set on the right does not have the $\delta_i\cdot \delta_i'=0$ condition. Unlike the family $\{\Phi_{s,t}\}_{s,t}$. Therefore $\{\Phi_{s,t}'\}_{s,t}$ should in principle be an easier family to prove dimension results for then $\{\Phi_{s,t}\}_{s,t}.$ 
\end{remark}

\end{document}